\documentclass[11pt]{article}

\usepackage{amsmath}
\usepackage{amssymb}
\usepackage{amscd}
\usepackage{amsthm}

\numberwithin{equation}{section}
\newtheorem{lemma}{Lemma}[section]
\newtheorem{consequence}{Corollary}[section]
\newtheorem{theorema}{Theorem}[section]

\newtheorem*{theor31}{Theorem \ref{Td>pi}}
\newtheorem*{theor32}{Theorem \ref{normsemiplane}}
\theoremstyle{definition}
\newtheorem{definition}{Definition}[section]
\newtheorem{counterexample}{Counterexample}[section]
\theoremstyle{remark}
\newtheorem{remark}{Remark}[section]
\newcommand{\N}{\mathbb N}

\newcommand{\R}{\mathbb R}

\newcommand{\E}{\mathbb E}

\newcommand{\Td}{\operatorname{Td}}
\newcommand{\Hd}{\operatorname{Hd}}
\newcommand{\dist}{\operatorname{dist}}
\newcommand{\const}{\operatorname{const}}
\newcommand{\Id}{\operatorname{Id}}

\begin{document}

\title{Tits geometry on ideal boundaries of Busemann non-positively curved space}
\author{P.D. Andreev}

\maketitle

\large

\section{Introduction}

In this paper, we consider Busemann non-positively curved spaces (shortly Busemann spaces). 
The term {\it Busemann space} was introduced by B. Bowditch in \cite{Bow}, the general 
geometric information on Busemann non-positively curved spaces can be found in \cite{Pa}. 
The class of Busemann spaces contains all $CAT(0)$-spaces and strongly convex normed spaces.

When $X$ is a complete locally compact $CAT(0)$-space, its geometry depends badly on the geometric boundary at infinity $\partial_\infty X$ and Tits metric $\Td$ on it.
Busemann's curvature non-positivity condition is weaker then Alexandrov's one. This leads to definite specialties of the geometry at infinity in Busemann spaces. Firstly, there are two natural approaches for the definition of the geometric boundary at infinity. In the $CAT(0)$-case the two approaches gives the same result, but when $X$ is Busemann space the results can be essentially different. We are to consider two different ideal boundaries --- horofunction (or metric) one and geodesic one. Secondly, even the two boundaries coincide, there are no natural way to define a metric on $\partial_\infty X$ with properties of Tits metric.

We propose the following trick that allows using the properties of Tits metric without the definition of the metric itself. Note that there are two key values of Tits metric on ideal boundary of $CAT(0)$-space. The values are $\pi$ and $\pi/2$. Conditions for ideal points $\xi, \eta \in \partial_\infty X$ under which inequalities
\begin{equation}\label{Tdgepi}
\Td(\xi, \eta) \ge \pi,
\end{equation}
\begin{equation}\label{Td<pi}
\Td(\xi, \eta) \le \pi
\end{equation}
and similar inequalities comparing Tits distance with $\pi/2$ hold, can be described geometrically without using Tits distance. In fact, mentioned inequalities can be considered as the collection of binary relations on the boundary $\partial_\infty X$. Consequently, one can define analogous collection of binary relations for Busemann space. 

Here we introduce a collection of binary relations of type \eqref{Tdgepi}, \eqref{Td<pi} etc. We prove that if $X$ is a proper Busemann space, some properties of Tits metric remain true for these binary relations. We use the notation of type $\Td(\xi, \eta) \le \pi$ to indicate that ideal points $\xi, \eta \in \partial_gX$ satisfy corresponding relation in analogy with $CAT(0)$-situation. One should not think that such notation means a comparison of metric function $\Td$ with $\pi$. The notation $\Td$ means only that the pair $(\xi, \eta)$ belongs to appropriate subset of $\partial_gX \times \partial_gX$.

The paper is structured as follows.

In Section \ref{prelim}, we recall some necessary facts from Busemann non-positively curved spaces theory. Also we specify definitions of horofunction and geodesic com\-pac\-ti\-fi\-ca\-tions $\overline{X}_h = X \cup \partial_hX$ and $\overline{X}_g = X \cup \partial_gX$. We establish relations between the compactifications.

In Section \ref{Titspi}, we introduce the collection of binary relations that generalize comparison the Tits distance with $\pi$ on the geodesic boundary of complete locally compact Busemann space. Here we prove the following two theorems generalizing known properties of Tits distance.

\begin{theor31}
Let $X$ be a proper Busemann space, and $\xi, \eta \in \partial_gX$ geodesic ideal points.
If $\Td(\xi, \eta) > \pi$, then there exists a geodesic $a : \R \to X$ with ends $a(-\infty) = \eta$ and $a(+\infty) = \xi$.
\end{theor31}

\begin{theor32}
Let $X$ be a proper Busemann space.
Given a geodesic $a: \R \to X$ with endpoints $\eta  = a(-\infty)$ and $\xi  = a(+\infty)$ passing throw $a(0)  = o$, the following conditions are equivalent.
\begin{enumerate}
\item $\Td(\xi, \eta) = \pi$; 
\item there exist horofunctions $\Phi$ centered in $\xi$ and $\Psi$ centered in $\eta$, such that the intersection of horoballs
\begin{equation}
\mathcal{HB}(\Phi, o) \cap \mathcal{HB}(\Psi, o)
\end{equation}
is unbounded;
\item there exists a normed semiplane in $X$ with boundary $a$.
\end{enumerate}
\end{theor32}

Here the horoball
$$\mathcal{HB}(\Phi, o) = \{y \in X\ |\ \Phi(y) \le \Phi(o)\}$$
is the sublevel set of the horofunction $\Phi$.

In Section \ref{Titspi2}, the collection of binary relations analogous to comparison of Tits metric with value $\pi/2$ is introduced. These relations are defined as subsets of $\partial_hX \times \partial_gX$. When $X$ is $CAT(0)$-space, the definition agrees with standard interpretation for the inequalities of type $\Td \le \pi/2$ etc. We also prove two versions of statement generalizing triangle inequality connecting relations $\Td(\xi, \eta) \le \pi$ with relations $\Td([\Phi], \theta) \le \pi/2$. Ambiguity of the formulation for such "triangle inequality"\ is a consequence of the relation $\Td([\Phi], \theta) \le \pi/2$ asymmetry. There is third possible way to formulate the "triangle inequality". Counterexample \ref{counter} shows that in general such third version of "triangle inequa\-li\-ty"\ is wrong.

In Section \ref{titstopology}, we apply the relation $\Td([\Phi], \theta) \le \pi/2$ to study the geometry of horoballs at infinity. If $\Phi$ is a horofunction then corresponding horoball at infinity $\mathcal{HB}_\infty(\Phi)$ can be presented as the set of points $\xi \in \partial_gX$ such that $\Td([\Phi], \xi) \le \pi/2$. We prove that $\mathcal{HB}_\infty(\Phi)$ is exactly the intersection of the boundary $\partial_gX$ with the closure $\overline{\mathcal{HB}(\Phi, y)}_g$ of arbitrary horoball $\mathcal{HB}(\Phi, y)$ in geodesic compactification $\partial_gX$. We also prove corresponding statement for horospheres at infinity in geodesically complete proper Busemann space $X$. In the case of horospheres the inclusion
$$
\mathcal{HS}_\infty(\Phi) \subset \overline{\mathcal{HS}(\Phi, y)}_g \setminus \mathcal{HS}(\Phi, y)
$$
can be exact.

\section{Busemann spaces and their ideal compactifications}\label{prelim}

In this section, we recall necessary basic facts from Busemann spaces theory and describe two  constructions of their boundary at infinity. We refer the reader to \cite{BBI}, \cite{BH} and \cite{Pa} for details in geometry of geodesic spaces and non-positively curved spaces.

\begin{definition}\label{defBuSp}
Let $(X, d)$ be a geodesic space. We use notation $|xy|$ for the distance $d(x, y)$ between points $x, y \in X$. The segment connecting points $x, y \in X$ will be denoted $[xy]$. The space $X$ is called \textit{Busemann space} (Busemann {\it non-positively curved} space) if for any two segments $[xy]$ and $[x' y']$ with corresponding affine parameterizations $\gamma: [a, b] \to X$, $\gamma': [a', b'] \to X$, the function $D_{\gamma, \gamma'}: [a,b]\times [a',b'] \to \R$  defined by
$$D_{\gamma, \gamma'}(t, t') = |\gamma(t)\gamma'(t')|$$
is convex. Equivalently, the geodesic space $X$ is Busemann space if for any three points $x, y, z \in X$, the midpoint $m$ between $x$ and $y$ and the midpoint $n$ between $x$ and $z$ satisfy the inequality
\begin{equation}\label{busem}
|mn| \le \frac12 |yz|.
\end{equation}
\end{definition}
Several more statements equivalent to the Definition \ref{defBuSp} are listed in \cite[Chapter 8]{Pa}. It follows easily from the definition \ref{defBuSp} that every Busemann space is contractible and every two points $x,y \in X$ are connected by the unique segment $[xy]$ in $X$. 

Every complete geodesic $c: \R \to X$ is embedding of the real line $\R$ to $X$. We call the image $c(\R)$ {\it straight line} in the space $X$. {\it The ray} is a geodesic $c: \R_+ \to X$, where $\R_+ = [0, +\infty)$. The rays $c, c' : \R_+ \to X$ are called {\it complement} if the map $a: \R \to X$ defined as $a(t) = c(t)$ for $t \ge 0$ and $a(t) = c'(-t)$ for $t \le 0$ represents a complete geodesic in $X$.

Simplest examples of Busemann spaces are $CAT(0)$-spaces and normed spaces with strongly convex norm.

\begin{definition}
The Haussdorff distance $\Hd(A, B)$ between closed subsets $A, B \subset X$ is
$$\Hd(A, B) = \inf\{\varepsilon\ |\ A \subset O_\varepsilon(B), B \subset O_\varepsilon (A)\}, $$
where
$$O_\varepsilon (C) = \{y \in X\ |\ \dist(y, C) < \varepsilon \}$$
denotes the $\varepsilon$-neighbourhood of the set $C\subset X$.

Straight lines $a,b \subset X$ are called \textit{parallel} if their Haussdorff distance is finite:
$$\operatorname{Hd}(a,b) < +\infty.$$
\textit{The normed strip} in the space $X$ is by definition the subset in $X$ isometric to the strip between two straight lines in normed plane. Every normed strip is bounded by two parallel straight lines forming \textit{the boundary} of the normed strip in $X$.
\end{definition}

\begin{lemma}[W. Rinow, \cite{Rin}, pp. 432, 463, \cite{Bow}, Lemma 1.1 and remarks]\label{lemmarinowa}
Every two parallel straight lines in Busemann space $X$ bound the normed strip in $X$.
\end{lemma}

Let $X$ be a proper Busemann space. Then it has two natural constructions of compactification. We call the first compactification $\overline{X}_g$ {\it geodesic} and the second one $\overline{X}_h$ {\it horofunction} compactification. When $X$ is $CAT(0)$-space, the two constructions give the same result in the following sense: the identity map $\Id: X \to X$ has continuation to a homeomorphism $\overline{X}_h \to \overline{X}_g$. In general Busemann space $X$ the compactifications $\overline{X}_g$ and $\overline{X}_h$ can be essentially different.

Now we specify explicit definitions.

\begin{definition}\label{geocomp}
Let $X$ be a non-compact proper Busemann space.
We say that geodesic rays  $c, d: \mathbb R_+ \to X$ in $X$ are \textit{asymptotic} if their Hausdorff distance is finite:
$$ \operatorname{Hd}(c(\mathbb R_+), d(\mathbb R_+)) < +\infty.$$
The asymptoticity relation $\sigma$ is equivalence on the family  $\mathcal{G}(\mathbb R_+, X)$ of geodesic rays in $X$. The factorset $\partial_g X = \mathcal{G}(\mathbb R_+, X)/\sigma$ forms the set of \textit{geodesic ideal points}. Given the ray $c: \mathbb R_+ \to X$ we denote $c(+\infty)$ corresponding geodesic ideal point. If $x \in X$ and $\xi \in \partial_gX$, then $[x \xi]$ denotes the ray from $x \in X$ in the class $\xi$. The ray $[x\xi]$ always exists and is unique.

We now define {\it cone topology} on the union $\overline{X}_g = X \cup \partial_g X$ as following. Fix a basepoint $o \in X$. By definition, the sequence $\{x_i\}_{i=1}^\infty \subset \overline{X}_g$ converges to the point $x \in \overline{X}_g$ in the sense of cone topology if
$$\lim_{i \to \infty}|ox_i| = |ox|$$
and the sequence
\begin{equation}\label{potoc}
c_i : [0, |ox_i|] \to X
\end{equation}
of naturally parameterized segments (rays) $[ox_i]$ converges to the naturally parameterized segment (ray) $[ox]$. Such a topology on the set $\overline{X}_g$ does not depend on the choice of the basepoint $o$. The space $\overline{X}_g$ with the cone topology is compact. It is called {\it geodesic compactification} of the space $X$. The set of ideal points $\partial_gX$ forms the {\it geodesic ideal boundary}.
\end{definition}

The cone topology restricted to the boundary $\partial_gX$ has the base of neighbourhoods of the point $\xi \in \partial_gX$ consisting of open sets
\begin{equation}\label{coneneib}
\mathcal U_{\xi, K}  = \{\eta \in \partial_gX \ |\ |c_\xi(K)c_\eta(K)| < 1\},
\end{equation}
which are defined for all $K > 0$. Here $c_\xi, c_\eta : \R_+ \to X$ are natural parameterizations of rays $c_\xi=[o\xi]$ and $ c_\eta= [o\eta]$.

From the other hand, given arbitrary non-compact proper metric space $(X, d)$, its horofunction (metric) compactification is defined as following.

\begin{definition}\label{horocomp}
Let $C(X, \R)$ be the space of continuous function with the topology of uniform convergence on bounded sets. Denote $C^\ast(X, \mathbb R)  = C(X, \mathbb R)/\{\operatorname{const}\}$ a factor-space of $C(X, \R)$ by the subspace of constants. Fix a basepoint $o \in X$ and identify every point $x \in X$ with corresponding distance function $d_x$: 
$$d_x(y) = d(x, y) - d(o, x).$$
The correspondence $x \to \{d_x + c\ |\ c = \const\}$ defines the embedding $\nu: X \to C^\ast(X, \mathbb R)$. The space $X$ is identified with its image $\nu(X) \in C^\ast(X, \R)$.  The {\it horofunction compactification} $\overline{X}_h$ is by definition the closure of the image $\nu(X)$ in $C^\ast(X, \mathbb R)$. The boundary $\partial_hX = \overline{X}_h \setminus X$ is called {\it horofunction boundary}, the functions forming the boundary $\partial_hX$ are called \textit{horofunctions}. Every ideal point of the horofunction boundary is a class of horofunctions which differ from each other by constants.
\end{definition}

This constructive approach to the horofunction compactification of the space $X$ is due to M. Gromov (\cite{Gro}). The embedding $\nu: X \to C^\ast(X, \R)$ was introduced also in \cite{Kur}[Ch. 2].

\begin{definition}
Given horofunction $\Phi$ on the space $X$, \textit{the horoball} corresponding to the point $y \in X$ is by definition the sublevel set
\begin{equation}\label{horobb}
\mathcal{HB}(\Phi, y) = \{x \in X\ |\ \Phi(x) \le \Phi(y)\}.
\end{equation}
The boundary of the horoball $\mathcal{HB}(\Phi, y)$ 
\begin{equation}\label{hoross}
\mathcal{HS} = \partial\mathcal{HB}(\Phi, y) = \{x \in X\ |\ \Phi(x) = \Phi(y)\}
\end{equation}
is called \textit{horosphere}.
\end{definition}

When $X$ is a proper Busemann space, the compactifications in definitions \ref{geocomp} and \ref{horocomp} satisfy the inequality $\overline{X}_g \le \overline{X}_h$ in the following sense. The identity map $\Id_X: X \to X$ has continuation to continuous surjection $\pi_{hg}: \overline{X}_h \to \overline{X}_g$. The map $\pi_{hg}$ can be non-injective (cf. \cite{An}, \cite{An2}). In particular, if the point $\xi \in \partial_gX$ is represented by the ray $c: \mathbb R_+ \to X$ such that $c(+\infty) = \xi$, then its pre-image $\pi_{hg}^{-1}(\xi)$ contains the class $[\beta_c] \in \partial_hX$ of the ray Busemann function $\beta_c$ defined for any $y \in X$ by the equality
$$\beta_c(y) = \lim_{t \to \infty} (|c(t)y| - t).$$

\begin{definition}
The point $\xi \in \partial_gX$ is called \textit{regular}, if $\pi_{hg}^{-1}(\xi)$ is one-point set
$$\pi_{hg}^{-1}(\xi) = \{[\beta_c]\}.$$
\end{definition}

It follows easily from the compactness of the space $\overline{X}_h$ and Hausdorff property of $\overline{X}_g$ that the map $\pi_{hg}$ is closed: the image of any closed subset in $\overline{X}_h$ is closed in $\overline{X}_g$. As a corollary, the map $\pi_{hg}$ satisfies the following "weak openness"\ property:

\begin{lemma}\label{weakopenness}
For any point $\xi \in \partial_gX$ and any neighbourhood $\mathcal U$ of its pre-image $\pi_{hg}^{-1}(\xi) \subset \partial_hX$ there exists a neighbourhood $\mathcal V$ of $\xi$ in $\partial_gX$ such that
$$\mathcal V \subset \pi_{hg}(\mathcal U).$$
\end{lemma}

\proof The image $\pi_{hg}(\partial_hX \setminus \mathcal U)$ is closed in $\partial_gX$, hence the set
$$\mathcal V = \partial_gX \setminus \pi_{hg}(\partial_hX \setminus \mathcal U)$$
is open. This set is desired neighbourhood of $\xi$. \qed

\begin{remark}
M. Rieffel in \cite{Ri1} introduces the notion of metric compactification for the proper non-compact metric space (see also \cite{WW}). This compactification is equivalent to the horofunction one described above. The term "horofunction compactification"\ was introduced by C. Walsh in \cite{Wa}.
\end{remark}

\section{Tits relations for the value $\pi$}\label{Titspi}

From now on the space $X$ will be a non-compact proper Busemann space.

Let rays $c, d: \R_+ \to X$ with common beginning
$$c(0) = d(0) = o$$ 
represent points $\xi  = c(+\infty)$ and $\eta  = d(+\infty)$ in the boundary $\partial_gX$. 
It follows from Busemann curvature non-positivity and triangle inequality that the function
$$\delta_{o, \xi, \eta}(t)  = \frac{|c(t)d(t)|}{2t}$$
is non-decreasing on $\R_+$ and bounded from above by $1$. Hence it has the limit
$$\delta_o(\xi, \eta)  =  \lim_{t \to +\infty}\delta_{o, \xi, \eta}(t)\le 1.$$

\begin{lemma}
Let rays $c, c' : \R_+ \to X$ be asymptotic in the direction $\xi \in \partial_gX$ and rays $d$ and $d': \R_+ \to X$ be asymptotic in the direction $\eta \in \partial_gX$. Let $c(0) = d(0) = o$ and $c' (0) = d' (0) = o'$. Then 
$$\delta_{o}(\xi, \eta) = \delta_{o'}(\xi, \eta).$$
\end{lemma}

\proof It follows from the metric convexity and asymptoticity of rays that
$$|c(t)c' (t)|, |d(t)d' (t)| \le |oo'|.$$
Hence the triangle inequality gives
$$|c(t)d(t)| - 2 |oo'| \le |c' (t)d' (t)| \le |c(t)d(t)| +2|oo'|.$$
For an arbitrary $\varepsilon > 0$ put
$$T = \frac{|oo'|}{\varepsilon}.$$
Then for any $t > T$ the inequality holds
$$\left|\frac{|c(t)d(t)|}{2t} - \frac{|c' (t)d' (t)|}{2t}\right| < \varepsilon.$$
This proves the claim of the lemma. \qed

\begin{remark}
When $X$ is $CAT(0)$-space, the inequality 
\begin{equation}\label{xieta<pi}
\delta_{o}(\xi, \eta) < 1
\end{equation}
holds iff the angle between ideal points $\xi$ and $\eta$ satisfies the inequality
$$\angle(\xi, \eta) < \pi.$$
In that case we have the equality for Tits distance
$$\Td(\xi, \eta) = \angle(\xi, \eta).$$
\end{remark}

The inequality \eqref{xieta<pi} has a characterization in terms of quasigeodesics. First recall the definition.

\begin{definition}\label{defquasigeodes}
Given numbers $a \ge 1$ and $b \ge 0$, the map $f: X \to Y$ between metric spaces $(X, d_X)$ and $(Y, d_Y)$ is called \textit{$(a, b)$-quasi-isometric} if for all $x, y \in X$ the two-sided inequality holds
\begin{equation}\label{quasiisometr}
\frac1a d_X(x, y) - b \le d_Y(f(x), f(y)) \le a d_X(x,y) + b.
\end{equation}
The subset $B \subset Y$ is \textit{$\varepsilon$-net} in $Y$ with $\varepsilon > 0$, if its $\varepsilon$-neighbourhood $N_\varepsilon(B)$ contains $Y$: $Y \subset \mathcal N_\varepsilon(B)$. 
The map $f: X \to Y$ is called \textit{$(a, b)$-quasi-isometry}, if it is $(a, b)$-quasi-isometric and the image $f(X)$ forms $\varepsilon$-net in $Y$ for some $\varepsilon> 0$. The map $f$ is called \textit{quasi-isometric} (correspondingly \textit{quasi-isometry}) if it is $(a, b)$-quasi-isometric (correspondingly $(a, b)$-quasi-isometry) for some $a \ge 1$ and $b \ge 0$. 
\textit{$(a, b)$-quasigeodesic} in the metric space $X$ is $(a, b)$-quasi-isometric map $f: I \to X$ for some $I \subset \R$. We also call \textit{$(a,b)$-quasigeodesic} the image  $f(I) \subset X$ in this map.
\end{definition}

Let two ideal points $\xi, \eta \in \partial_gX$ in the proper Busemann space $X$ with basepoint $o \in X$ be given. Consider rays $c  = [o\xi]$ and $d  = [o, \eta]$ with corresponding natural parameterizations $c, d: \R_+ \to X$.
We define the map $f: \R \to X$ by the equality
\begin{equation}\label{fcd}
f(t)  =
	\left\{
		\begin{array}{ll}
			c(t) & \mbox{ when } t \ge 0\\
			d(-t) & \mbox{ when } t \le 0
		\end{array}			
	\right..
\end{equation}

\begin{lemma}
Given ideal points $\xi, \eta \in \partial_gX$, the equality $\delta_o(\xi, \eta) = \pi$ holds iff for any $\varepsilon > 0$ there exists a number $b \ge 0$ such that the map \eqref{fcd} is $(1 +\varepsilon, b)$-quasigeodesic in $X$.
\end{lemma}

\proof 
\textit{Necessity.}
Let $\delta_o(\xi, \eta) = \pi$. Since the rays $c$ and $d$ are geodesic, we need only to verify the condition \eqref{quasiisometr} for arbitrary numbers $-s, t \in \R, -s < 0 < t$, that is for points $d(s)$ and $c(t)$.

From the one hand, the triangle inequality gives $|d(s)c(t)| \le t+s$, and therefore the right inequality in \eqref{quasiisometr} holds. From the other hand, the condition $\delta_{o}(\xi, \eta) = 1$ means that for any $\varepsilon > 0$ there exists a number $T > 0$, such that
$$|c(t)d(t)| > \frac{2t}{1+\varepsilon_1}$$
for all $t > T$, where $\varepsilon_1$ satisfies to the equation
$$\frac{1}{1+ \varepsilon} = \frac{1 - \varepsilon_1}{1 + \varepsilon_1}.$$
We have for $t > T$ and $s \le t$ 
\begin{align*}
|c(t)d(s)| &\ge |c(t)d(t)| - t + s > \frac{2t}{1+\varepsilon_1} - \frac{(t - s)(1 + \varepsilon_1)}{1+\varepsilon_1} = \\
&= \frac{t+s}{1+\varepsilon_1} - \frac{\varepsilon_1(t-s)}{1+\varepsilon_1} \ge (t+s)\frac{1 - \varepsilon_1}{1 + \varepsilon_1} = \frac{t+s}{1 + \varepsilon}.
\end{align*}
Similarly,
$$|c(t)d(s)| \ge \frac{t+s}{1 + \varepsilon}$$
for $s > T$ and $t \le s$.

Since the segment $[-T, T]$ is compact, there exists a number $b_1 \ge 0$ such that, if $T \ge t \ge s$, then
$$|c(s)d(t)| > \frac{s+t}{1+\varepsilon} - b_1.$$
Thus, inequality
\begin{equation}\label{b1}
|c(s)d(t)| > \frac{s+t}{1+\varepsilon} - b_1
\end{equation}
holds for all $t \ge s$.
Similar arguments show, that there exists a number $b_2$ such that  
\begin{equation}\label{b2}
|c(s)d(t)| > \frac{s+t}{1+\varepsilon} - b_2
\end{equation}
for all $s \ge t$. Take $b  = \max\{b_1, b_2\}$. Then we receive from inequalities \eqref{b1} and \eqref{b2} demanded
$$|c(s)d(t)| > \frac{s+t}{1+\varepsilon} - b.$$

\textit{Sufficiency.} Let $\delta_o(\xi, \eta) < \pi$ and $\varepsilon > 0$ be such that
$$\frac{1}{1+\varepsilon} > \delta_{o, \xi, \eta}.$$
We claim that there is no $b \ge 0$ for which the map $f$ is $(1 + \varepsilon, b)$-quasigeodesic.

Indeed, 
$$|c(t)d(t)| \le 2\delta_{o,\xi, \eta}t$$
for all $t \ge 0$. Given $b \ge 0$ there exists $T > 0$ such that
$$b < 2\left(\frac{1}{1+\varepsilon} - \delta_{o, \xi, \eta}\right)T.$$
Hence for all $t \ge T$
$$|f(-t)f(t)| \le \frac{2t}{1 + \varepsilon} - b,$$
contradicting to the left-side inequality in the definition of $(1 + \varepsilon, b)$-quasigeodesic.
\qed

\begin{remark}
The statement of the lemma does not mean that the map \eqref{fcd} is $(1, b)$-quasigeodesic for some $b > 0$. The simplest counterexample can be constructed as following. Fix a number $\alpha$, \  $1 < \alpha < 2$. Take the Euclidean plane $\E^2$ with coordinates $(x, y)$. The space $X$ is  constructed by deleting from $\E^2$ two convex domains bounded by curves $x^2 = \pm |y|^\alpha$ and gluing two Euclidean half-planes to their places.  Consider the following curve $\gamma: \R \to X$ in $X$. The image $\gamma(\R)$ consists of two semiparabolas  $x^2 = |y|^\alpha$ in the half-plane $x > 0$. The parameterization of $\gamma$ is natural. The curve $\gamma$ satisfies the conditions of the lemma (it connects two ideal points with angle distance $\pi$). But the computation shows that there is no $b > 0$ such that $\gamma$ is $(1, b)$-quasigeodesic. To see this, fix $b \ge 0$ and consider the function 
$$f(x) = x^{2-\alpha}(x^{2\frac{\alpha-1}{\alpha}} - 2b).$$
It is increasing to $+\infty$ when $x^{2\frac{\alpha-1}{\alpha}} > 2b$. Hence there is $x_0 > 0$ satisfying
$$f(x_0) > b^2.$$
Denote $y_0 = x_0^{2/\alpha}$. Then
$$x_0^2 > 2by_0 + b^2$$
and
$$\sqrt{x_0^2 + y_0^2} > y_0 + b.$$
The points $A(x_0, -y_0)$ and $B(x_0, y_0)$ belong to $\gamma$. The distance between them is $2y_0$. At the same time, the length of the path $\gamma$ from $A$ to $B$ is greater then
$$2|OA| = 2\sqrt{x_0^2 + y_0^2} > 2y_0 + b.$$
Hence $\gamma$ is not $(1,b)$-quasigeodesic.
\end{remark}

Next, we apply the notation $\Td$ to define five binary relations on $\partial_gX$. The notation $\Td(\xi, \eta)$ denotes Tits distance between ideal points $\xi, \eta\in \partial_\infty X$ in $CAT(0)$-spaces theory. In our case the notation of type $\Td(\xi, \eta) < \pi$ is related only for a binary relation, not for any metric.

\begin{definition}\label{deftits}
Let $(X, o)$ be a pointed proper Busemann space.
Given arbitrary ideal points $\xi, \eta \in \partial_gX$ we define the following binary relations:
\begin{itemize}
\item $\Td(\xi, \eta) < \pi$ if $\delta_o(\xi, \eta) < \pi$
\item $\Td(\xi, \eta) \le \pi$, if for any neighbourhoods $U(\xi)$ and $V(\eta)$ of this points in the sense of cone topology on $\partial_gX$ there exist points $\xi'  \in U(\xi)$ and $\eta'  \in V(\eta)$ with $\Td(\xi' , \eta' ) < \pi$;
\item $\Td(\xi, \eta) \ge \pi$, if $\Td(\xi, \eta) < \pi$ does not hold;
\item $\Td(\xi, \eta) > \pi$ if $\Td(\xi, \eta) \le \pi$ does not hold;
\item $\Td(\xi, \eta)  = \pi$ if $\Td(\xi, \eta) \ge \pi$ and $\Td(\xi, \eta) \le \pi$ hold simultaneously.
\end{itemize}
\end{definition}

The definition \ref{deftits} is correlated with the standard definition of Tits metric when $X$ is $CAT(0)$. Obviously, $\Td(\xi, \eta) < \pi \Rightarrow \Td(\xi, \eta) \le \pi$ and $\Td(\xi, \eta) > \pi \Rightarrow \Td(\xi, \eta) \ge \pi$. It follows directly from the definition that 
$$ \Td(c(-\infty), c(+\infty) \ge \pi$$
for any geodesic $c: \R \to X$.

Moreover, we have the following lower semicontinuity property of $\Td$-relation with respect to the cone topology. If a consequence $\{\xi_n\} \subset \partial_gX$ converges in the sense of cone topology to the point  $\xi \in \partial_gX$, a consequence  $\{\eta_n\} \subset \partial_gX$ converges to the point $\eta \in \partial_gX$, and if for all $n \in \N$ we have $\Td(\xi_n, \eta_n) \le \pi$, then $\Td(\xi, \eta) \le \pi$.

Now we study some properties of introduced relations.

\begin{lemma}\label{interpoint}
Suppose that $\Td(\xi, \eta) < \pi$, $c  = [o\xi]$ and $c'   = [o \eta]$. Let $m_t$ be an arbitrary point of the segment $[c(t)c' (t)]$ where $t > 0$. Then $|om_t| \to \infty$ when $t \to \infty$.
\end{lemma}

\proof Suppose that there exists a consequence $t_n \to \infty$ for which distances $|om_{t_n}|$ are bounded. Choose an accumulation point $m \in X$ for the family $m_{t_n}$. Such accumulation point does exist because the space $X$ is proper. The natural parameterizations of segments $[m_{t_n}c(t_n)]$ and $[m_{t_n}c' (t_n)]$ converge correspondingly to natural parameterizations of rays $[m\xi]$ and $[m\eta]$ uniformly on bounded domains in $\R$. Since each pair of mentioned segments have opposite directions in the point $m_{t_n}$, the rays $[m\xi]$ and $[m\eta]$ are complement. A contradiction. \qed

\begin{lemma} \label{caphoroballs}
Let $c, c': \R \to X$ be rays with $c(0) = c' (0) = o$. Let the intersection of horoballs $\mathcal{HB}(\beta_c, o) \cap \mathcal{HB}(\beta_{c'}, o)$ is compact. Then
$$\Td(c(+\infty), c' (+\infty) \ge \pi$$ 
and there exists a geodesic $a: \R \to X$ with $a(+\infty) = c(+\infty)$ and $a(-\infty) = c' (+\infty)$.
\end{lemma}

\proof Denote
$$T  = \inf\{t \in \R_+\ |\ \mathcal{HB}(\beta_c, c(t)) \cap \mathcal{HB}(\beta_{c'}, c' (t)) \ne \varnothing\}.$$
It follows from the horofunctions continuity and compactness of the intersection 
$$\mathcal{HB}(\beta_c, o) \cap \mathcal{HB}(\beta_{c'}, o)$$
that the set 
$$ M(T)  = \mathcal{HB}(\beta_c, c(T)) \cap \mathcal{HB}(\beta_{c'}, c' (t))$$
is non-empty. It contains minimum points of the function $\max\{\beta_c(x), \beta_{c'}(x)\}$ in the intersection $\mathcal{HB}(\beta_c, o) \cap \mathcal{HB}(\beta_{c'}, o)$.
Choose an arbitrary point $m \in M(T)$. Then rays $[mc(+\infty)]$ and $[mc' (+\infty)$ are complement. \qed

{\sloppy

\begin{remark}
The intersection of horoballs $\mathcal{HB}(\beta_c, o) \cap \mathcal{HB}(\beta_{c'}, o)$ for complement rays $c, c': \R_+ \to X$ may be compact even under the condition $\Td(c(+\infty),
c'(+\infty)) = \pi$ (cf. example in \cite{An2}). If $X$ is $CAT(0)$, such situation is impossible.
\end{remark}

}

The two following theorems generalize statements 1 and 3 of Proposition 9.11 in \cite{BH}.

\begin{theorema}\label{Td>pi}
Let $X$ be a proper Busemann space.
If $\Td(\xi, \eta) > \pi$, then there exists a geodesic $a : \R \to X$ with ends $a(-\infty) = \xi$ and $a(+\infty) = \eta$.
\end{theorema}

\proof Suppose that the following two statements hold simultaneously:\par
{\bf A}) $\Td(\xi, \eta) > \pi$\\
 and \par
{\bf B}) the straight line in $X$ with endpoints $\xi$ and $\eta$ does not exist.

Fix a basepoint $o \in X$, rays $c  = [o\xi]$ and $c'  = [o\eta]$ with natural parameterizations $c, c': \R_+ \to X$ and arbitrarily large number $K > 0$. The condition  {\bf A}) allows to choose the number $K$ such that
$$|c(K+t)c'(K)| > K + t + 1$$
 for all $t \ge 0$, and consequently
\begin{equation}\label{betac>0}
\beta_c(x)> 0
\end{equation}
for all $x \in B(c'(K), 1)$.
Analogously we may assume that
\begin{equation}\label{betac'>0}
\beta_{c'}(y) > 0
\end{equation}
for all  $y \in B(c(K), 1)$. Here $\beta_c$ and $\beta_{c'}$ are Busemann functions defined within rays $c$ and $c'$ correspondingly. In such a choice of $K$ we have \eqref{betac>0} for all $x \in B(c'(\alpha K), \alpha)$ and $\eqref{betac>0}$ for all $y \in B(c(\alpha K), \alpha)$, where $\alpha \ge 1$ is arbitrary.

By the curvature non-positivity, there exists a number $T > K$ such that for any $t \ge T$ and points  $a, b\in X$ with $|ac(t)|, |bc'(t)| \le 2$, the following holds. If $\gamma_a : [0, |oa|] \to X$ is natural parameterization of the segment $[x_0 a]$ and $\gamma_b : [0, |x_0b|] \to X$ is natural parameterization of the segment $[x_0 b]$, then 
\begin{equation}\label{cagamma}
|c(K)\gamma_a(K)| < 1
\end{equation} 
and
\begin{equation}\label{cbgamma}
|c' (K)\gamma_b(K)| < 1.
\end{equation}

Take points $a'_1, b'_1 \in [c(T)c' (T)]$ on distances 1 from corresponding endpoints $c(T)$ and $c' (T)$:
$$|c(T)a' _1| = |c' (T)b' _1| = 1.$$

\begin{figure}
	\begin{center}
	\begin{picture}(170,200)
		\qbezier(10,10)(30,50)(50,120)
		\qbezier(50,120)(60,150)(60,210)
		\qbezier(10,10)(40,50)(80,90)
		\qbezier(80,90)(130,140)(170,140)
		\qbezier(57,150)(80,120)(116,120)
		\put(10,10){\circle*{3}}
		\put(57,150){\circle*{3}}
		\put(116,120){\circle*{3}}
		\put(67,139){\circle*{3}}
		\put(98,122){\circle*{3}}
		\put(24,145){$c(T)$}
		\put(120,112){$c'(T)$}
		\put(70,142){$a'_1$}
		\put(97,128){$b'_1$}
		\put(55,220){$\xi$}
		\put(180,135){$\eta$}
		\put(0,0){$o$}
		\end{picture}
	\end{center}
\caption{}\label{premidpoi}
\end{figure}

It may happen that one of the two points $a'_1$ or $b'_1$ belongs to the ray $[o\xi]$ or to the ray $[o\eta]$ correspondingly. But both equalities\\
$$a'_1 = c(T-1)$$
 and
$$b'_1 = c'(T-1)$$
can not hold simultaneously, because in that case the union of the segment $[a'_1b'_1]$ with rays $[a'_1\xi]$ and $[b'_1\eta]$ is a straight line with endpoints $\xi$ and $\eta$. If $|oa'_1| > |o b'_1|$, then $a'_1 \ne c(T-1)$. So we can assume (after renotation if necessary) that  $|oa'_1| \ge |ob'_1|$ and $a'_1 \notin [o\xi]$. Denote $a_1  = a'_1$ and $b_1$ the point of the segment $[b'_1c'(T)]$, on the distance $\tau_1  = |x_0a_1|$ from $o$.

Next, for all natural $n \ge 2$ we define points $a_n$ and $b_n$ and numbers $\tau_n$. For this, consider the following functions
$$\phi_n,\phi'_n: [0, |c(nT)c'(nT)|] \to \R_+.$$ 
The values $\phi_n(s)$ and $\phi'_n(s)$ for $s \in [0, |c(nT)c'(nT)|]$ are defined as following. Let $y_{n,s}\in [c(nT)c'(nT)]$ be the point on the distance $s$ from $c(nT)$, and $z_{n,s}\in [oy_{n,s}]$ be on the distance $T$ from $o$. Then $\phi_n(s)$ and $\phi'_n(s)$ are by definition
$$\phi_n(s)  = |c(T)z_{n,s}|,$$
and
$$\phi'_n(s)  = |c'(T)z_{n,s}|.$$

\begin{figure}
	\begin{center}
	\begin{picture}(170,200)
		\qbezier(5,5)(15,25)(25,60)
		\qbezier(25,60)(30,75)(30,105)
		\qbezier(30,105)(32,150)(29,210)
		\qbezier(5,5)(20,25)(40,45)
		\qbezier(40,45)(65,70)(85,70)
		\qbezier(85,70)(150,75)(170,72)
		\qbezier(31,180)(50,120)(160,73)
		\qbezier(5,5)(35,65)(67,129)
		\put(5,5){\circle*{3}}
		\put(31,180){\circle*{3}}
		\put(160,73){\circle*{3}}
		\put(67,129){\circle*{3}}
		\put(36,65){\circle*{3}}
		\put(28,70){\circle*{3}}
		\put(57,59){\circle*{3}}
		\put(-5,185){$c(nT)$}
		\put(0,70){$c(T)$}
		\put(140,60){$c'(nT)$}
		\put(57,50){$c'(T)$}
		\put(65,135){$y_{n,s}$}
		\put(40,65){$z_{n,s}$}
		\put(28,220){$\xi$}
		\put(180,70){$\eta$}
		\put(-5,-5){$o$}
	\end{picture}
	\end{center}
\caption{}
\end{figure}

The function $\phi_n(s)$ is continuous, $\phi_n(0)=0$ and $\phi_n(|c(nT)c' (nT)|)= |c(T)c' (T)|$. Similarly, the function $\phi' _n(s)$ is continuous,  $\phi'_n(0)=|c(T)c' (T)|$  and $\phi'_n(|c(nT)c' (nT)|)=0$. Denote  $a' _n \in [c(nT)c' (nT)]$ the farest point from $c(nT)$ such that $\phi_n(s) \le 1$ for all $s \le |c(nT)a' |$, and $b' _n$ the farest point from $c' (nT)$ such that $\phi' _n(s) \le 1$ for all $s \ge |c(nT)b' _n|$.
Set 
$$\tau_n  = \max\{|oa' _n|, |ob' _n|\}.$$ 
If $|oa' _n| = \tau_n \ge |ob' _n|$, then we redenote $a_n  = a' _n$ and denote $b_n$ the point of the segment $[b' _nc' (nT)]$ with $|ob_n|= \tau_n$. In that case the point $a_n$ does not belong to the ray $[o\xi]$ and 
\begin{equation}\label{phi=1}
\phi_n(|c(nT)a_n|) = 1.
\end{equation}

From the other hand, if $|oa'_n| < |ob'_n| = \tau_n$, then $a_n$ is the point of the segment  $[c(nT)a'_n]$ on the distance $\tau_n$ from $o$ and we redenote $b_n  = b'_n$. In that case the point $b_n$ does not belong to the ray $[o\eta]$ and
\begin{equation}\label{phi'=1}
\phi'_n(|c(nT)b_n|) = 1.
\end{equation} 

With such a choice, the points of the segments $[oa_n]$ and $[ob_n]$ on the distance $T$ from $o$ are separated from points $c(T)$ and $c'(T)$ correspondingly on the distance not greater then $1$. Moreover, one of the two following possibilities holds necessarily. The consequence $a_n$ contains an infinitely many points satisfying the equality  \eqref{phi=1}, or the consequence $b_n$ contains infinitely many points satisfying the equality \eqref{phi'=1}. For definiteness, assume the first case.

Notice that the values of Busemann functions $\beta_c$ and $\beta_{c'}$ are negative in points $a_n$ and $b_n$ correspondingly:
$$\beta_c(a_n), \beta_{c'}(b_n) < 0,$$
however the values $\beta_{c'}$ and $\beta_c$ are positive (cf. \eqref{betac>0} and \eqref{betac'>0}):
$$\beta_{c'}(a_n), \beta_c(b_n) > 0.$$

Therefore, for each segment $[a_nb_n]$ one can find a point $m_n$ where $\beta_c(m_n) = \beta_c'(m_n)$. It follows from the estimation
$$|c(nT)c'(nT)| < 2nT = -\beta_c(c(nT)) - \beta_{c'}(c'(nT))$$
that
$$\beta_c(m_n) = \beta_{c'}(m_n) < 0.$$
Consequently, the point $m_n$ belongs to the intersection of the interiors of horoballs:
\begin{equation}\label{mnincap}
m_n \in \operatorname{Int}(\mathcal{HB}(\beta_c, o)) \cap \operatorname{Int}(\mathcal{HB}(\beta_{c'}, o)).
\end{equation}

Independently on the choice of the consequence $m_n$ with condition \eqref{mnincap}, any such a consequence $\{m_n\}_{n=1}^\infty$ has an accumulation point $m$ in the compact topological space $\overline{X}_g$. There are two possibilities:
\par 1) there exists a consequence $\{m_n\}_{n=1}^\infty$ satisfying the condition \eqref{mnincap} with accumulation point in $X$: $m \in X$,\\
or \par
2) each accumulation point $m$ for any consequence of type $\{m_n\}_{n=1}^\infty$ is infinite:  $m \in \partial_gX$. In this case all accumulation points for all consequences $\{m_n\}_{n=1}^\infty$ where $m_n \in [c(nT)c'(nT)]$ are infinite as well. In particular, all accumulation points for consequences $\{a_n\}_{n=1}^\infty$ and $\{b_n\}_{n=1}^\infty$ are infinite. 

The first possibility contradicts to the condition {\bf B}). The proof of this fact is analogous to that of Lemma \ref{interpoint}. 

We claim that  the second possibility contradicts to the condition {\bf A}).

To prove that, pick out a consequence $n_k$ such that a subconsequence of points $a_{n_{k}}$ converges in the cone topology to $\theta \in \partial_gX$ and for all $k$ the equality holds
$$\phi_{n_{k}} (a_{n_{k}}) = 1,$$
and subconsequence of points $b_{n_{k}}$ converges in the cone topology to $\zeta \in \partial_gX$. By the construction we have $\theta \in \mathcal U_{\xi, K}$ and $\zeta \in \mathcal U_{\eta, K}$. We will show that rays $p  = [o\theta]$ and $q  = [o\zeta]$ satisfy to the condition
$$\lim_{t \to \infty}\frac{|p(t)q(t)|}{2t} < 1.$$

For each $n \in \N$ the ball $B(a_n, |a_nm_n|)$ is contained in the horoball $\mathcal{HB}(\beta_c, o)$, and the ball $B(b_n, |b_nm_n|)$ in the horoball $\mathcal{HB}(\beta_{c'}, o)$. Hence
$$|a_nb_n| \le -\beta_c(a_n) - \beta_{c'}(b_n).$$

Consequently, for each $t \le \tau_n$ the following holds. If $y_t \in [oa_n]$, and $z_t \in [ob_n]$ are points on the distance $t$ from $o$, then
$$|y_tz_t| \le -\beta_c(y_t) - \beta_{c'}(z_t).$$
When $k \to \infty$, we have in the limit
$$|p(t)q(t)| \le -\beta_c(p(t)) - \beta_{c'}(q(t)).$$

We estimate the value $- \beta_{c'}(q(t))$ by
\begin{equation}\label{-beta'}
-\beta_{c'}(q(t)) \le t.
\end{equation}

Next, estimate the value $-\beta_{c}(p(t))$ using the equality \eqref{phi=1}.
Put 
$$z  = p(T) = \lim\limits_{k \to \infty}z_{n_{k}, |c(n_kT)a_{n_{k}}|}.$$

For each $k \in \N$ the point $z_{n_{k}, |c(n_kT)a_{n_{k}}|}$ belongs to the compact intersection of spheres $S(o, T) \cap S(c(T), 1)$. Consequently
$$z \in S(o, T) \cap S(c(T), 1).$$
The function $-\beta_c$ when restricted to $S(o, T) \cap S(c(T), 1)$ attains its maximum in the point $w \in S(o, T) \cap S(c(T), 1)$. The maximum value is
\begin{align*}
-\beta_c(w) = \max\{-\beta_c(v)\ |\ v \in S(o, T) \cap S(c(T), 1)\} < \\
<\max\{-\beta_c(v) \ |\ v \in S(o, T)\} = T
\end{align*}
and
$$-\beta_c(z) < T.$$
Since the function $\beta_c$ is convex, the inequality holds
\begin{equation}\label{-betac}
-\beta_c(q(t)) < \frac{-\beta_c(z)}{T}\cdot t
\end{equation}
 for all $t > T$. Combination of the estimations \eqref{-beta'} and \eqref{-betac} gives the inequality
$$\frac{|p(t)q(t)|}{2t} \le \frac{\left(1 + \frac{-\beta_c(z)}{T}\right)}{2} <1$$
 for all $t > T$. Hence
$$\Td(\theta, \zeta) < \pi.$$
Since the number $K$ was chosen arbitrarily, it follows 
$$\Td(\xi, \eta) \le \pi.$$
A contradiction with the condition {\bf A}). \qed

\begin{consequence}
Let $X$ be a proper Busemann space. If points $\xi, \eta \in \partial_gX$ satisfy to the relation $\Td(\xi, \eta) > \pi$, then there exist their neighbourhoods $\mathcal U(\xi), \mathcal U(\eta) \subset \partial_gX$, such that any $\theta \in \mathcal U(\xi)$ and $\zeta \in \mathcal U(\eta)$ are connected by a straight line $a: \R \to X$ with endpoints $a(+\infty) = \theta$ and $a(-\infty) = \zeta$.
\end{consequence}

\proof

By definition of the relation $\Td(\xi, \eta) > \pi$ there exist neighbourhoods $\mathcal U(\xi), \mathcal U(\eta) \subset \partial_gX$, such that if $\theta \in \mathcal U(\xi)$ and $\zeta \in \mathcal U(\eta)$, then $\Td(\theta, \zeta) > \pi$. Hence the claim for points $\theta$ and $\zeta$ follows from proven Theorem \ref{Td>pi}. \qed

\begin{lemma}\label{noquasi}
Let $\alpha$ be strictly convex normed plane and $\xi, \eta \in \partial_g\alpha$ be not opposite endpoints of any straight line in $\alpha$. Then there does not exist $(1, d)$-quasigeodesic $b: \R \to X$ with $b(t) \to \xi$ when $t \to +\infty$, $b(t) \to \eta$ when $t \to -\infty$. Here $d \ge 0$ is arbitrary, convergence is meant in the sense of the cone topology in $\overline{\alpha}_g$.
\end{lemma}

\proof Suppose in contrary that the mentioned quasigeodesic does exist. Without lost of generality we may assume that the quasigeodesic $b: \R \to \alpha$ is continuous. In fact, given $(1, d)$-quasigeodesic one can obtain continuous $(1, d')$-quasigeodesic for some $d' \ge d$ by replacing some its pieces by straight segments with affine parameterizations.

Pick a point $o  = b(0)$. It follows from the continuity of the quasigeodesic $b$ that for any $t > 0$ there exists $\tau \in \R$ with $|ob(\tau)| = t$. The definition of $(1, d)$-quasigeodesic gives the inequality
$$|t - \tau| \le d.$$
Consider points $p \in [o\eta]$ and $q \in [o\xi]$ on distances $|op| = |oq| = 1$ from $o$. The distance $|pq|$ is
\begin{equation}\label{pq2delta}
|pq| = 2 - \Delta < 2
\end{equation}
for some  $\Delta > 0$.
From the other hand, let the numbers $\tau_n < 0$ and  $\tau'_n > 0$ be such that $|ob(\tau_n)| = |ob(\tau'_n)| = n$. Set points $p_n \in [ob(\tau_n)]$ and $q \in [o\tau'_n]$ on distances $|op_n| = |oq_n| = 1$ from $o$. Then
$$|p_nq_n| = \frac1n |b(\tau_n)b(\tau'_n)| \ge \frac1n(\tau'_n - \tau_n - d) \ge 2 - \frac{3d}{n}.$$
For $n > 6d/\Delta$ the inequality holds $|p_nq_n| \ge 2 - \frac12\Delta$. Combination with \eqref{pq2delta} contradicts to the convergence conditions $b(t) \to \xi$ when $t \to +\infty$ and $b(t) \to \eta$ when $t \to -\infty$.\qed

\begin{definition}
\textit{The normed semiplane} in the metric space $X$ is by definition the subset in $X$ isometric to a half-plane in the two-dimensional normed space. When $X$ is Busemann space, each its normed semiplane is convex subset isometric to a half-plane in normed space with strongly convex norm. Given isometry $i: \bar \alpha \to \bar \alpha' \subset X$, where $\bar \alpha$ is a half-plane bounded by straight line $a$ in the normed space $V^2$, we say that the image $i(a)$ bounds the normed semiplane $\bar \alpha'$ in $X$. Obviously, $i(a)$ is straight line in $X$.
\end{definition}

Let a geodesic $a$ bounds a normed semiplane in $X$. Then it follows from Lemma\ref{noquasi} that $\Td(a(+\infty), a(-\infty)) = \pi$. The following theorem gives more complicated statement.

\begin{theorema}\label{normsemiplane}
Let $X$ be a proper Busemann space. Given a geodesic $a: \R \to X$ with endpoints $\xi  = a(+\infty)$ and $\eta  = a(-\infty)$ passing throw $a(0)  = o$, the following conditions are equivalent.
\begin{enumerate}
\item $\Td(\xi, \eta) = \pi$; 
\item there exist horofunctions $\Phi$ centered in $\xi$ and $\Psi$ centered in $\eta$, for which the intersection of horoballs
\begin{equation}\label{intersechoroballs}
\mathcal{HB}(\Phi, o) \cap \mathcal{HB}(\Psi, o)
\end{equation}
is unbounded;
\item  $a$ bounds a normed semiplane in $X$.
\end{enumerate}
\end{theorema}

\proof
$1 \Rightarrow 2$.
Suppose that for any two horofunctions $\Phi$ and $\Psi$ centered in $\xi$ and $\eta$ correspondingly, the intersection \eqref{intersechoroballs} is bounded. The pre-images $\pi_{hg}^{-1}(\xi)$ and $\pi_{hg}^{-1}(\eta)$ are closed in $\partial_hX$, and consequently compact. We claim that all intersections \eqref{intersechoroballs} are uniformly bounded in  $X$: they are contained in some ball $B(o, R)$.

Indeed, suppose that there exist consequences 
$$\{[\Phi_n]\}_{n=1}^\infty \subset \pi_{hg}^{-1}(\xi)$$ 
and 
$$\{[\Psi_n\}_{n=1}^\infty \subset \pi_{hg}^{-1}(\eta),$$ 
such that
$$\Phi_n(o) = \Psi_n(o) = 0$$
and
$$R_n = \min\{R |\mathcal{HB}(\Phi_n, o) \cap \mathcal{HB}(\Psi_n, o) \subset B(o, R)\} \to \infty.$$
Using the compactness of the boundary $\partial_hX$ we may assume (passing to a subconsequences, if necessary), that the consequences of horofunctions $\Phi_n$ and $\Psi_n$ converge to horofunctions $\Phi$ and $\Psi$ correspondingly. They also satisfy conditions $\pi_{hg}([\Phi]) = \xi$ and $\pi_{hg}(\Psi) = \eta$.
Let the consequence $\{a_n\}_{n=1}^\infty\subset X$ satisfies to equality $|oa_n| = R_n$ and
$$\Phi_n(a_n), \Psi_n(a_n) \le 0.$$
Then 
$$\Phi_n(x), \Psi_n(x) \le 0.$$
for all $x \in [oa_n]$. We may assume that $\{a_n\}_{n=1}^\infty$ converges to the infinite point $\theta \in \partial_gX$. In that case the natural parameterizations of segments $[oa_n]$ converge to the natural parameterization of the ray $[o\theta]$: every point $x\in [o\theta]$ is the limit point of the consequence $x_n \in [oa_n]$ with $|ox_n| = |ox|$. Since $\Phi_n(x_n), \Psi_n(x_n) \le 0$, then for any $\varepsilon > 0$ there exists a natural $N$, such that the inequalities
$$\Phi_n(x), \Psi_n(x) \le \varepsilon$$
hold for all $n \ge N$. We have in the limit $\Phi(x), \Psi(x) \le 0$. Since $x$ is an arbitrary point of the ray $[o\theta]$, we obtain the contradiction with the boundedness of the intersection \eqref{intersechoroballs}.

Let the number $R > 0$ be such that for any horofunction $\Phi$ in the class $[\Phi] \in \pi_{hg}^{-1}(\xi)$ and  $\Psi$ in the class $[\Psi] \in \pi_{hg}^{-1}(\eta)$ the inclusion holds
$$\mathcal{HB}(\Phi, o) \cap \mathcal{HB}(\Psi, o) \subset B(o, R).$$
For every $[\Phi] \in \partial_hX$ we consider those a horofunction  $\Phi$ for which $\Phi(o) = 0$. In such a choice of horofunctions $\Phi$ and $\Psi$ for every $x \in S(o, 2R)$ the inequality holds
$$\Phi(x) + \Psi(x) > 0.$$
The continuous function $\Phi + \Psi$ attains its maximum in the sphere $S(o, 2R)$ and the maximal value is positive. Hence there are neighbourhoods $O_\Psi(\Phi)$ and $O_\Phi(\Psi)$ in the space $C(X, \R)$ such that $f(x) + g(x) > 0$ for any $f \in O_\Psi(\Phi)$ and $g \in O_\Phi(\Psi)$ and for all $x \in S(o, 2R)$. Denote $\mathcal U_\Psi(\Phi) = p(O_\Psi(\Phi)) \cap \partial_hX$ and $\mathcal U_\Phi(\Psi) = p(O_\Phi(\Psi)) \cap \partial_hX$, where $p: C(X, \R) \to C_*(X, \R)$ is a projection map.

Fix $[\Phi] \in \partial_hX$ and consider a covering of the compact set $\pi_{gh}^{-1}(\eta) \subset \partial_hX$
by open sets of type $\mathcal U_\Phi(\Psi)$ for all $[\Psi] \in \pi_{hg}^{-1}(\eta)$. Pick a finite subcovering:
$$\pi_{gh}^{-1}(\eta) \subset \bigcup_{i = 1}^n \mathcal U_{\Phi_i}(\Psi_i).$$
Denote
$$\mathcal U_\eta(\Phi) = \bigcup_{i = 1}^n \mathcal U_{\Psi_i}(\Phi_i)$$
and
$$\mathcal U_\Phi(\pi_{hg}^{-1}(\eta)) = \bigcup_{i = 1}^n \mathcal U_{\Phi_i}(\Psi_i).$$
Now consider the covering of the compact set $\pi_{hg}^{-1}(\xi)$ by open sets of type $\mathcal U_\eta(\Phi)$ for all $[\Phi] \in \partial_hX$. Taking its finite subcovering we get an open set
$$\mathcal V_+ = \bigcup_{j = 1}^N \mathcal U_\eta(\Phi_j) \supset \pi_{hg}^{-1}(\xi)$$
and open set
$$\mathcal V_- = \bigcap_{j = 1}^N \mathcal U_{\Phi_j}(\eta),$$
with the following property. If $[\Theta_+] \in \mathcal V_+$ and $[\Theta_-] \in \mathcal V_-$, then for all $x \in S(o, 2R)$
$$\Theta_+(x) + \Theta_-(x) > 0.$$
Consequently,
\begin{equation}\label{twoTheta}
\mathcal{HB}(\Theta_+, o) \cap \mathcal{HB}(\Theta_-, o) \subset B(o, 2R).
\end{equation}
From Lemma \ref{weakopenness}, there exists neighbourhoods $\mathcal U_+ \subset \partial_gX$ of the point $\xi$ and $\mathcal U_- \subset \partial_gX$ of the point $\eta$ with the following inclusions
$$\mathcal U_+ \subset \pi_{hg}(\mathcal V_+)$$ 
and
$$\mathcal U_- \subset \pi_{hg}(\mathcal V_-).$$
Since the map $\pi_{hg}$ is continuous, we may assume that 
$\mathcal V_\pm = \pi_{hg}^{-1}(\mathcal U_\pm)$.
The horofunction $\Theta_+$ attains its minimum in the compact set \eqref{twoTheta} in some point $y_0$. The point $y_0$ belongs to the boundary of the set \eqref{twoTheta}, hence $y_0 \in \mathcal{HS}(\Theta_-, o)$ and 
$$\Theta_-(y_0) = \Theta_-(o) = \max \Theta_-|_{\mathcal{HB}(\Theta_+, o) \cap \mathcal{HB}(\Theta_-, o)}.$$
Moreover,  $y_0$ is the minimum point for the function $\Theta_+$ on $\mathcal{HS}(\Theta_-, o)$, since $\Theta_+(y_0) \le \Theta_+(o)$ и $\Theta_+(q) > \Theta_+(o)$ for all $q \in \mathcal{HB}(\Theta_-, o) \setminus \mathcal{HBH}(\Theta_+, o)$.
Put $\theta_\pm  = \pi_{hg}([\Theta_\pm])\in \partial_gX$. Consider rays $[y_0\theta_-]$ and $[y_0\theta_+]$. 
For any point $p \in [y_0\theta_+]$ we have the inequality $\Theta_-(p) \le \Theta_-(y_0) + |py_0|$. 
From the other hand, the projection of the point $p$ (the nearest to $p$ point) to the closed convex set $\mathcal{HB}(\Theta_-, o)$  is $y_0$. Indeed, if $z \in \mathcal{HB}(\Theta_-, o)$, then
$$|pz| \ge \Theta_+(z) - \Theta_+(p) \ge \Theta_+(y_0) - \Theta_+(p) = |py_0|,$$
and since any point in $X$ has unique projection to the bounded compact convex subset, hence the point $y_0$ is namely the projection of $p$. 

Take a point $q \in [y_0\theta_-]$. Denote $m$ the intersection point of the segment $[pq]$ with the horosphere $\mathcal{HS}(\Theta_-, o)$. We have
$$|pq| = |pm| + |mq| \ge |py_0| -\Theta_-(o) = |py_0| + |y_0 q| \ge |pq|.$$
Now we conclude that $m = y_0$ and rays $[y_0\theta_+]$ and $[y_0 \theta_-]$
complement each other up to a complete geodesic connecting  $\theta_-$ with $\theta_+$. 
This contradicts to the relation $\Td(\xi, \eta) \le \pi$.

$2 \Rightarrow 3$.

Assume that the intersection \eqref{intersechoroballs} is not compact. Then its closure in $\overline{X}_g$ contains infinite points: 
$$\overline{\mathcal{HB}(\Phi, o) \cap \mathcal{HB}(\Psi, o)}_g \cap \partial_gX \ne \varnothing.$$
Each ray $[o\theta]$ with
$$\theta \in \overline{\mathcal{HB}(\Phi, o) \cap \mathcal{HB}(\Psi, o)}_g \cap \partial_gX,$$
is contained in $\overline{\mathcal{HB}(\Phi, o) \cap \mathcal{HB}(\Psi, o)}_g \cap \partial_gX$.
Moreover, given a point $x \in X$, the intersection $\mathcal{HB}(\Phi, x) \cap \mathcal{HB}(\Psi, x)$ is non-compact as well, and
$$\overline{\mathcal{HB}(\Phi, x) \cap \mathcal{HB}(\Psi, x)}_g \cap \partial_gX = \overline{\mathcal{HB}(\Phi, o) \cap \mathcal{HB}(\Psi, o)}_g \cap \partial_gX.$$

We assume that $\Phi(o) = \Psi(o) =0$.
Given $K > 0$ denote $N_K(a)$ the $K$-neighborhood of the geodesic $a$. The function $|\Phi + \Psi|$ is bounded from above by $2K$ on  $\overline{N_K(a)}$. We claim that the following statements hold.
\begin{enumerate}
\item\label{2.01.09 12.10} The restriction of $\Phi + \Psi$ to $\overline{N_K(a)}$ attains its minimum in some point $b_K$;
\item for any minimum point  $b_K \in \overline{N_K(a)}$ of $\Phi+\Psi$ the rays $[b_K\xi]$ and $[b_K \eta]$ are complement;
\item the function $\Phi + \Psi$ is constant on the rays $[b_K\xi]$ and $[b_K \eta]$;
\item \label{2.01.09 12.11} the point $b_K$ can be chosen in $\mathcal{HB}(\Phi, o) \cap \mathcal{HB}(\Psi, o) \cap S(o, K)$ on the distance $K$ from $a$.
\end{enumerate}

Prove them.
Given arbitrary $K' > K$, the function $\Phi+ \Psi$ attains its minimum in the compact subset $B(o, K') \cap \overline{N}_K(a)$ in some point $p_{K'}$. At least one of the two rays $[p_{K'}a(+\infty)]$ or $[p_{K'} a(-\infty)]$ intersects the interior of the ball $B(o, K')$ in the point different from $p_{K'}$. For definiteness, assume that the beginning part of the ray $c =  [p_K a(+\infty)$ passes in the interior of $B(o,K')$. Then $\Phi + \Psi$ is constant in $[p_{K'} a(+\infty)$, because $(\Phi + \Psi)(c(t))$ is convex, bounded from above and non-decreasing at the initial segment $[0, t_0]$, which corresponds to the part of the ray in $B(o,K')$.  Similarly, $\Phi + \Psi$ is constant along the ray $[c(t_0) a(-\infty)]$. Consequently, $p_{K'} \in [c(t_0) a(-\infty)]$, and the rays $[p_{K'}a(-\infty)]$ and $[p_{K'}a(+\infty)]$ are complement. The union $a'_{p_{K'}} = [p_{K'} a(-\infty)] \cup [p_{K'}a(+\infty)]$ is the complete geodesic parallel to $a$.  

Since the sum $\Phi + \Psi$ is non-increasing on the ray $[p_{K'} \theta]$, we can assume that $\dist(a, a'_{p_{K'}}) = K$. The geodesic $a'_{p_{K'}}$ contains the point
$$b_K \in a'_{p_{K'}} \cap \mathcal{HB}(\Phi, o) \cap \mathcal{HB}(\Psi, o)$$
and
$$(\Phi+\Psi)(b_K) = (\Phi+\Psi)(p_{K'}) = \min\limits_{x \in B(o, K' \cap \overline{N}_K(a)}(\Phi + \Psi)(x)).$$

For each $K > 0$ the geodesics  $a$ and $a'_K = [b_K a(+\infty)] \cup [b_K a(-\infty)]$ bound the normed strip $\mathcal F_K$. It follows that the minimum of the value $\Phi + \Psi$ in $B(o, K') \cap \overline{N}_K(a)$ is the same for all $K' > K$. Hence, the point $b_K$ can be chosen the same for all for all $K' > K$. This proves all statements \ref{2.01.09 12.10} -- \ref{2.01.09 12.11} listed above.

Note that the function $\Phi + \Psi$ is linear on each segment $[ob_K]$ for all $K > 0$. Hence, if
$$b_K, b'_K \in S(o, K) \cap \overline{N}_K(a) \cap (\mathcal{HB}(\Phi, o) \cap \mathcal{HB}(\Psi, o),$$
and $\nu_K, \nu'_K : [0, K] \to X$ are the natural parameterizations of the segments $[ob_K]$ and $[ob'_K]$ correspondingly, then
$$(\Phi+\Psi)(\nu_K(t)) = (\Phi + \Psi)(\nu'_K(t)$$
for all $t \in [0, K]$.
Let $\theta\in \partial_X$ be the accumulation point for the family $b_K$ when $K \to +\infty$. Then $$(\Phi + \Psi)(c(t)) = (\Phi + \Psi)(\nu_K(t))$$
for all $K > 0$ and $t \in [0, K]$. Here $c: \R_+ \to X$ (correspondingly, $\nu_K : [0,K] \to X$) is the natural parameterization of the ray $[o \theta]$ (correspondingly, segment $[ob_K]$). By this reason we may assume that  $b_K = c(K)$ for all $K > 0$. Under such assumption, the normed strips $\mathcal F_K$ are ordered by inclusion: $\mathcal F_{K_1}\subset \mathcal F_{K_{2}}$, if $K_1 \le K_2$. The union
$$\bar \alpha = \bigcup_{K > 0} \mathcal F_K$$
is the required normed semiplane.

$3 \Rightarrow 1$.
By Definition \ref{deftits}, $\Td(\xi, \eta) \ge \pi$. Suppose that $\Td(\xi, \eta) > \pi$. Consider a neighbourhood $\mathcal U_+$ of the point $\xi$ and a neighbourhood $\mathcal U_-$ of the point $\eta$ in $\partial_gX$ such that any pair of points in these neighbourhoods are the pair of endpoints for some geodesic in $X$. Draw rays $[o \theta_+]$ and $[o \theta_-]$ in directions of some ideal points $\theta_\pm \in \mathcal U_\pm$ correspondingly differ from $\xi$ and $\eta$ in the normed semiplane with boundary $a$. Let $b: \R \to X$ be a geodesic in $X$ with endpoints $\theta_+$ and $\theta_-$. The projection $p \circ b$ of the geodesic $b$ to the normed semiplane $\bar \alpha$ bounded by $a$ represents $(1, d)$-quasigeodesic, where $d = 2\max \operatorname{dist}(b, \bar \alpha)$. 
In fact, the projection $p$ to the semiplane is a submetry: $|p(x)p(y)| \le |xy|$ for all $x, y \in X$. Hence
$$|(p \circ b)(s) (p\circ b)(t)| \le |s-t|$$
for all $s, t \in \R$. From the other hand,
$$|s-t| = |b(s) b(t)| \le |b(s) (p \circ b)(s)| + |(p\circ b)(s)(p \circ b)(t)| + |(p \circ b)(t) b(t)| \le$$
$$\le |(p\circ b)(s)(p \circ b)(t)| + 2\max \operatorname{dist}(b, \bar \alpha),$$
and consequently
$$|(p\circ b)(s)(p \circ b)(t)| \ge |s-t| - d.$$
Notice that the map $p \circ b$ represents $(1, d)$-quasigeodesic within the interior metric of the semiplane $\bar \alpha$.
But Lemma \ref{noquasi} states that the normed semiplane with strongly convex norm admits no $(1, d)$-quasigeodesic with endpoints different from endpoints of its boundary. A contradiction.
\qed

\section{Tits relations for the value $\pi/2$}\label{Titspi2}

Here we introduce similar collection of binary relations corresponding to the angle value $\pi/2$. The following lemma serves as motivation for the Definition \ref{pi2} below.

\begin{lemma}
Let $X$ be a proper $CAT(0)$-space, $\xi, \theta \in \partial_\infty X$. Given arbitrary points $x, y \in X$ consider rays $c = [x\xi]$ and $d = [y \theta]$ with natural parameterizations $c, d: \R_+ \to X$ correspondingly. Let $\beta_c: X \to \R$ be a Busemann function defined by the ray $c$.
\begin{enumerate}
\item The inequality $\Td(\xi, \eta) \le \pi/2$ holds iff the function $\beta_c \circ d: \R_+ \to \R$ is non-increasing on $\R_+$.
\item The inequality $\Td(\xi, \eta) < \pi/2$ holds iff the function $\beta_c \circ d: \R_+ \to \R$ decreases sublinearly, that is if there exist numbers $k < 0$ and $b \in \R$, such that for all $t \in \R_+$ the inequality holds
$$(\beta_c \circ d)(t) \le kt+b.$$
\end{enumerate}
\end{lemma}

\proof Firstly, consider the case when $x = y$. In that case the inequality $\Td(\xi, \eta) \le \pi/2$ is equivalent to the condition 
\begin{equation}\label{piphagor}
|c(s)d(t)|^2 \le s^2 + t^2
\end{equation}
for all $s, t \ge 0$. Let be $\Td(\xi, \eta) \le \pi/2$. Fix a number $t \ge 0$ and start to enlarge the value $s$ unboundedly. Then the inequality \eqref{piphagor} implies
$$\beta_c(d(t)) = \lim\limits_{s \to +\infty} (|c(s)d(t)| - s) \le
\lim\limits_{s \to +\infty}(\sqrt{s^2+t^2} - s) = 0.$$
Hence the function $\beta_c \circ d$ is non-positive. It follows from its convexity, that $\beta_c \circ d$ is non-increasing function.

Conversely, if the function $\beta_c \circ d$ is non-increasing, then $\beta_c(d(t)) \le 0$ for all $t \ge 0$. Fix  $t_0, s_0 \ge 0$. Then for any $\varepsilon > 0$ there is $S > s_0$ such that for all $s > S$ the inequality holds\
$$|c(s)d(t)| < s + \varepsilon.$$
If there exists $s_1$ such that $|c(s_1)d(t_0)| \le s_1$, then we obtain from $CAT(0)$-inequality for the triangle $\bigtriangleup (xc(s_1)d(t_0))$
$$|c(s_0)d(t_0)| < s_0^2 + t_0^2.$$
If we can not find suitable value $s_1$, then for arbitrary $s_1 > S$ denote $z$ the point of the segment $[c(s_1)d(t_0)|$ with $|zc(s_1)| = s_1$. Now it follows from  $CAT(0)$-inequality for the triangle $\bigtriangleup (xc(s_1)z)$ and the triangle inequality that
$$|c(s_0)d(t_0)| < s_0^2 + t_0^2 + \varepsilon.$$
Since the choice of $\varepsilon > 0$ for values $s_0, t_0$ was arbitrary, we have
$$|c(s_0)d(t_0)| < s_0^2 + t_0^2.$$
The values $s_0$ and $t_0$ was also chosen arbitrarily, consequently we have $\Td(\xi, \eta) \le \pi/2$ for points $\xi, \eta \in \partial_\infty X$.

The inequality $\Td(\xi, \eta) < \pi/2$ is equivalent to the condition that there exists a number $\lambda > 0$, with condition
\begin{equation}\label{cosineq}
\frac{s^2+ t^2 - |c(s)d(t)|^2}{2st} > \lambda
\end{equation}
for all $s, t > 0$. The inequality \eqref{cosineq} is equivalent to 
$$|c(s) d(t)| < \sqrt{s^2 + t^2 - 2 \lambda st},$$
and we obtain for fixed $t > 0$ that
$$\beta_c(d(t)) = \lim_{s \to \infty} (|c(s) d(t)| - s) \le -\lambda t.$$ 
Finally, in the case $x=y$ the strong inequality $\Td(\xi, \eta) < \pi/2$ is equivalent to the following one
$$|c(s) d(t)| < kt, $$
for all $s, t > 0$, where the number $k$ can be taken as $k = -\lambda/2$.

If the points $x$ and $y$ does not coincide, consider the ray $d'$ with beginning $x$ asymptotic to $d$: $d' = [x\theta]$. It is clear that the function $\beta_c \circ d$ is non-increasing (correspondingly sublinearly decreasing) iff the function $\beta_c \circ d'$ is such. This proves the lemma in the general case. \qed

\begin{lemma}\label{mono}
Let $X$ be a proper Busemann space and $\Phi : X \to \R$~--- horofunction.
\begin{enumerate}
\item If for some ray $c: \R_+ \to X$ the function $\Phi \circ c$ is non-increasing, then for any ray $c' : \R_+ \to X$ asymptotic to $c$ the function $\Phi \circ c'$ is non-increasing as well.
\item If for some ray $c: \R_+ \to X$ the function $\Phi \circ c$ decreases sublinearly, then for any ray $c' : \R_+ \to X$ asymptotic to $c$ the function $\Phi \circ c'$ decreases sublinearly on $\R_+$ as well.
\end{enumerate}
\end{lemma}

\proof Notice that if the function $f: \R_+ \to \R$ is convex, it is non-increasing iff it is bounded from above by the value $f(0)$. The convexity of the horofunction $\Phi$ means that its restriction to any geodesic segment is convex. If the ray $c'$ is asymptotic to $c$, then there are estimations
$$|c(t)c' (t)| \le |c(0) c' (0)|$$
and 
\begin{equation}\label{PhiPhi'}
\Phi(c'(t)) \le \Phi(c(0)) + |c(0)c'(0)|.
\end{equation}
Hence, if $\Phi\circ c$ is non-increasing, then $\Phi \circ c'$ is bounded from above. Obviously, the maximal value is $(\Phi \circ c)(0)$. Consequently,  $\Phi \circ c'$ is also non-increasing.

Let $\Phi \circ c$ decreases sublinearly: 
$$\Phi(c(t)) < kt + b$$
for some $k < 0$ and $b \in \R$ and for all $t \ge 0$.
Then it follows from \eqref{PhiPhi'} that 
$$\Phi(c'(t) < kt + b + |c(0)c'(0)|,$$
for all $t \ge 0$ and consequently the claim of the lemma is true.
\qed

Now we are ready to define the new collection of binary relations. Every one of them is formally a subset in $\partial_hX \times \partial_gX$. We use the notation $\Td$ again. The notation reflects the analogy with the Tits metric.

\begin{definition}\label{pi2}
Let $X$ be a proper Busemann space. Given ideal points $[\Phi] \in \partial_hX$ and $\xi \in \partial_gX$ we define the following binary relations.
\begin{enumerate}
\item $\Td([\Phi], \xi) \le \pi/2$ if the horofunction $\Phi$ is non-increasing on some ray $c: \R_+ \to X$ with the endpoint $c(+\infty) = \xi$. In this case by the Lemma \ref{mono} any horofunction $\Phi'  \in [\Phi]$ is non-increasing on any ray $c' : \R_+ \to X$ asymptotic to $c$. 
\item $\Td([\Phi], \xi) < \pi/2$ if $\Phi$ decreases sublinearly on any ray $c: \R_+ \to X$ with endpoint $c(+\infty) = \xi$.
\item $\Td([\Phi], \xi) > \pi/2$ if not $\Td([\Phi], \xi) \le \pi/2$;
\item $\Td([\Phi], \xi) \ge \pi/2$ if not $\Td([\Phi], \xi) < \pi/2$;
\item $\Td([\Phi], \xi) = \pi/2$ if $\Td([\Phi], \xi) \ge \pi/2$ and $\Td([\Phi], \xi) \le \pi/2$ simultaneously.
\end{enumerate}
\end{definition}

Geometrically the condition $\Td([\Phi], \xi) > \pi/2$ means that rays in $\xi$-direction leave any horoball $\Phi \le \operatorname{const}$. The statement of the following lemma generalizes to the condition $\Td(\xi, \eta) \le \pi/2$ in Busemann spaces the property of low-semicontinuity with respect to the cone topology, known in $CAT(0)$-spaces case.

\begin{lemma}\label{limmm}
Let the sequence of horofunctions $\Psi_n$ converges (uniformly on bounded subsets) to the horofunction $\Phi$ and sequence of geodesic ideal points $\zeta_n$ converges in the cone topology to the point $\xi \in \partial_g X$. Let $\Td([\Psi_n], \zeta_n) \le \pi/2$ for all natural $n$. Then $\Td([\Phi], \xi) \le \pi/2$.
\end{lemma}

\proof
Fix a basepoint $o \in X$. Let $c_n, c: \R_+ \to X$ be natural parameterizations of rays $[o\zeta_n]$ and $[o\xi]$ correspondingly. Let numbers $t, \varepsilon > 0$ be arbitrary and  $N \in \N$ be such that for all $n > N$ the following conditions hold.
\begin{enumerate}
\item $|\Phi(c(t)) - \Psi_n(c(t))| < \frac{\varepsilon}{2}$ 
\\and
\item $|c(t)c_n(t)| < \frac{\varepsilon}{2}$.
\end{enumerate}
Then
$$\Phi(c(t)) < \Psi_n(c_n(t)) + \varepsilon \le \Psi_n(o) + \varepsilon = \Phi(o) + \varepsilon, $$
because the horofunction $\Psi_n$ is non-increasing on the ray $c_n$ for every natural $n$. Therefore, since $\varepsilon > 0$ is arbitrary,
$$\Phi(c(t)) \le \Phi(c(0)).$$
Since $t > 0$ is also arbitrary, the horofunction $\Phi$ is bounded: $\Phi(c(t)) \le \Phi(c(0))$ for all $t \in \R_+$. Consequently, $\Phi$ is non-increasing function on the ray $c$ and $$\Td([\Phi], \xi) \le \frac{\pi}{2}.$$
\qed

In cases when geodesic and metric compactifications of the space $X$ coincide, or when the ideal point $\eta \in \partial_gX$ is regular, we will write $\Td(\eta, \xi) > \pi/2$ and in the same way another relations from the Definition \ref{pi2} for the class of Busemann function $[\beta_\eta]\in \partial_hX$ which projects to $\eta$. Notice that the relations are not symmetric in general: it is possible to be true $\Td(\eta, \xi) > \pi/2$ and $\Td(\xi, \eta) < \pi/2$ simultaneously. For example, such pairs of points can be found in every non Euclidean normed space. Unfortunately, the author does not know, whether the relation $\Td(\xi, \eta) \le \pi/2$ implies $\Td(\xi, \eta) < \pi$. But there are two versions of "triangle inequality" for introduces relations. We formulate them in two following theorems.

\begin{theorema}\label{trnglemma1}
Let $\xi, \eta \in \partial_g X$ and $[\Phi] \in \partial_hX$ be ideal points such that
\begin{equation}\label{phixi}
\Td([\Phi], \xi) \le \pi/2
\end{equation} 
and
\begin{equation}\label{phieta}
\Td([\Phi], \eta) \le \pi/2.
\end{equation} 
Then $\Td(\xi, \eta) \le \pi$.
\end{theorema}

\proof If points $\xi$ and $\eta$ are not endpoints for a geodesic in $X$, the claim of the theorem is true by Theorem \ref{Td>pi}.

Suppose that there exists a geodesic $a: \R \to X$ with endpoints $a(-\infty)  = \eta$ and $a(+ \infty)  = \xi$. The function $\Phi \circ a$ is non-increasing and non-decreasing convex function, therefore it is a constant. We assume $\Phi|_a = 0$. 
Set $o = a(0)$. Denote $\theta = \pi_{hg}([\Phi])$ and consider a ray $c = [o\theta]$ with natural parameterization $c: \R_+ \to X$. 

For arbitrary $K > 0$ consider the ray $d'_K = [c(K)\xi]$ with natural parameterization $d'_K: \R_+ \to X$. The function $\Phi\circ d'_K$ is non-increasing on $\R_+$ because of the condition \eqref{phixi}. But since $|a(t)d'_K(t)| \le |a(0) d'_K(0)| = K$ for all $t > 0$, hence 
$$\Phi(d'_K(t)) \ge \Phi(a(t)) - K = \Phi(c(K)).$$ 
We conclude that both the function $\Phi$ and the distance function $|a(t)d'_K(t)|$ are constant on the ray $[c(K)\xi]$: 
$$\Phi(d'_K(t)) = -|a(t)d'_K(t)| = -K.$$
Analogously, the function $\Phi$ and the distance function $|a(-t)d''_K(t)|$ are constant on the ray $d''_K=[c(K)\eta]$ as well: 
$$\Phi(d''_K(t)) = -|a(-t)d''_K(t)| = -K.$$

We show that the rays $d'_K$ and $d_K''$ complement each other to the complete geodesic $d_K$ parallel to $a$. For this, consider points $d'_K(T)$ and $d''_K(T)$, where $T > 0$ is arbitrary and the midpoint $m$ of the segment $[d'_K(T)d''_K(T)]$. We have the following two estimations for the distance $|om|$. At first, it follows from the metric convexity in the space $X$
$$|om| = |a(0)m| \le \frac12(|a(T) d'_K(T)| + |a(-T)d''_K(-T)| = K.$$
From the other hand, since the point $c(K)$ is the projection of $o$ to the horoball $\mathcal{HB}(\Phi, c(K))$, then
$$|om| \ge |oc(K)| = K.$$
Consequently $|om| = K$, $m = c(K)$, and therefore the union  $d_K$ of rays $d_K'$ and $d_K''$ is a complete geodesic. Since all its points are on the same distance $K$ from $a$, we obtain that $a$ and $d_K$ are parallel. Denote $F_K$ the normed strip between $a$ and $d_K$. It is clear that $F_{K_{1}} \subset F_{K_{2}}$ when $K_1 < K_2$. Therefore, the union
$$\bar \alpha = \bigcup_{K > 0}F_K$$
is a normed semiplane in $X$ with boundary $a$. The relation $\Td(\xi, \eta) = \pi$ follows now from the Theorem \ref{normsemiplane}.

\qed

\begin{theorema}\label{triang2}
Let the ideal points $[\Phi], [\Psi] \in \partial_hX$ and the ideal point $\zeta \in \partial_gX$ be such that 
\begin{equation}\label{xizeta}
\Td([\Phi], \zeta) \le \pi/2
\end{equation}
and
\begin{equation}\label{etazeta}
\Td([\Psi], \zeta) \le \pi/2.
\end{equation}
Let $\xi  = \pi_{hg}([\Phi])$ and $\eta  = \pi_{hg}([\Psi])$. 
Then $\Td(\xi, \eta) \le \pi$.
\end{theorema}

\proof 
Suppose that endpoints of some geodesic $a:\R\to X$ are $a(+\infty) = \xi$ and $a(-\infty) = \eta$. Denote $a(0)  = o$ and $b = [o\zeta]$. By the condition, restrictions of functions $\Phi$ and $\Psi$ to $b$ are non-increasing. Hence the intersection of horoballs $\mathcal{HB}(\Phi, o) \cap \mathcal{HB}(\Psi, o)$ is non-compact. Applying Theorem \ref{normsemiplane}, we get that $a$ bounds a normed semiplane in $X$ and $\Td(\xi, \eta) = \pi$.

If the geodesic $a$ does not exist, the inequality $\Td(\xi, \eta) \le \pi$ is the corollary of the Theorem \ref{Td>pi}. \qed

A priori one can formulate another version for the triangle inequality: 
$$\Td(\xi, \zeta) \le \pi/2, \Td(\zeta, \eta) \le \pi/2 \Rightarrow \Td(\xi, \eta) \le \pi?$$
The following counterexample shows that the third version of the triangle inequality is not correct.

\begin{counterexample}\label{counter}

Consider three items $\alpha_1$, $\alpha_2$ and $\alpha_3$ of normed semiplane with coordinates $(x_1, x_2)$, $x_2 \ge 0$ and  norm
$\|(x_1, x_2)\|  = \sqrt[4]{x_1^4 + x_2^4}$

Glue them to the metric space with interior metric in the following way: the positive boundary ray $x_1$ of the semiplane $\alpha_1$ is glued to the negative boundary ray $-x_1$ of the semiplane  $\alpha_2$, the positive boundary ray $x_1$ of the semiplane $\alpha_2$ is glued to the negative boundary ray $-x_1$ of the semiplane $\alpha_3$ and finally, the positive boundary ray $x_1$ of the semiplane $\alpha_3$ is glued to the negative boundary ray $-x_1$ of the semiplane $\alpha_1$. The resulting space $X$ is Busemann space. Its geodesic and horofunction compactifications coincide: the surjection $\pi_{hg}: \overline{X}_h \to \overline{X}_g$ is a homeomorphism. There are ideal points $\xi, \eta, \zeta \in \partial_gX$, such that $\Td([\beta_\xi], \zeta) < \pi/2, \Td([\beta_\zeta], \eta) < \pi/2$, but $\Td(\xi, \eta) > \pi$. For example, such points are infinite points of the following rays: the point $\xi$ on the ray directed by the vector $(\cos\frac{5\pi}{6}, \sin\frac{5\pi}{6})$ in the semiplane $\alpha_3$, the point $\zeta$ on the ray with directing vector $(\cos\frac{\pi}{3}, \sin\frac{\pi}{3})$ in the semiplane $\alpha_1$ and the point $\eta$ on the ray with directing vector $(\cos(\frac{5\pi}{6}+\varepsilon), \sin(\frac{5\pi}{6} + \varepsilon))$, where
$$0 < \varepsilon < 
\arctan
\left(3^{\frac34}\right) - \frac{\pi}{3}
$$
in the semiplane $\alpha_1$.

\end{counterexample}

\section{Horoballs at infinity}\label{titstopology}

\begin{definition}\label{orisharnabeskti}
Let $\Phi:X \to \R$ be a horofunction that generates an ideal point $[\Phi] \in \partial_hX$. \textit{The horoball at infinity} with \textit{center} $[\Phi]$ is by definition a set
$$\mathcal{HB}_\infty(\Phi)  = \left\{\xi \in \partial_gX\ |\ \Td([\Phi], \xi) \le \frac{\pi}{2} \right\}.$$
Correspondingly, the set
$$\mathfrak{hb}_\infty(\Phi)  =  \left\{\xi \in \partial_gX\ |\ \Td([\Phi], \xi) < \frac{\pi}{2}\right\}$$
is called \textit{open horoball at infinity}, and the set
$$\mathcal{HS}_\infty(\Phi)  =  \left\{\xi \in \partial_gX\ |\ \Td([\Phi], \xi) = \frac{\pi}{2}\right\}$$
is \textit{horosphere at infinity} with center $[\Phi]$. 
\end{definition}

\begin{remark}
In general the set $\mathfrak{hb}_\infty(\Phi)$ is not the interior for $\mathcal{HB}_\infty(\Phi)$ in the sense of the cone topology and $\mathcal{HS}_\infty(\Phi)$ is not its boundary. These statements are false even with respect to Tits metric when $X$ is  $CAT(0)$-space. The reason of such effect is: the closed ball in Tits metric can be a component of linear connection. In this case it will be open set. For example, consider the space $X$ obtained by gluing by the boundary $\ell$ of the Euclidean semiplane $\bar \alpha$ and Lobachevskii semiplane $\bar \beta$. For infinite point $\xi \in \partial_\infty X$ on the ray perpendicular to $\ell$ in the Euclidean semiplane, the horoball at infinity $\mathcal{HB}_\infty(\beta_\xi)$ coincides with the boundary $\partial_\infty(\bar \alpha)$, and the horosphere at infinity $\mathcal{HS}_\infty(\beta_\xi)$ with $\partial_\infty (\ell)$. The horoball $\mathcal{HB}_\infty(\beta_\xi) = \partial_\infty(\bar \alpha)$ is open set in Tits metric.
\end{remark}

Now we give another description for horoballs and horospheres at infinity in the proper Busemann space $X$. Every horoball $\mathcal{HB}(\Phi, y)$ is a sublevel set \eqref{horobb} in $X$, and the horosphere $\mathcal{HS}(\Phi, y)$ is corresponding level set \eqref{hoross}.
When $X$ is contained in $\overline{X}_g$ with cone topology, all subsets in $X$ have their closures in $\overline{X}_g$.

\begin{lemma}\label{hb=hb}
\begin{enumerate}
\item Let $X$ be a proper Busemann space. Then
\begin{equation}\label{hotoball_infty}
\overline{\mathcal{HB}(\Phi, y)}_g = \mathcal{HB}(\Phi, y) \cup \mathcal{HB}_\infty(\Phi)
\end{equation}
\item
Let $X$ be a geodesically complete proper Busemann space. Then
\begin{equation}\label{horosphere_infty}
\mathcal{HS}_\infty(\Phi) \subset \overline{\mathcal{HS}(\Phi, y)}_g \setminus \mathcal{HS}(\Phi, y)
\end{equation}
\end{enumerate}
\end{lemma}

\proof
\begin{enumerate}
\item
From the definition of horoball at infinity, the inclusion $\xi \in \mathcal{HB}_\infty([\Phi])$ holds iff the horofunction $\Phi$ is non-increasing at the ray $[y\xi]$. Equivalently, $[y\xi] \subset \mathcal{HB}(\Phi, y)$. Consequently the equality \eqref{hotoball_infty}.

\item 
It follows from the definition of the horosphere at infinity that the inclusion $\xi \in \mathcal{HS}_\infty([\Phi])$ holds iff the horofunction  $\Phi$ is non-increasing and does not decreases sublinearly on the ray $[y\xi]$. In that case the ray $[y\xi]$ belongs to the horoball $\mathcal{HB}(\Phi, y)$ but may pass strongly in the interior of the horoball $\mathcal{HB}(\Phi, y)$ and be not a subset of $\mathcal{HS}(\Phi, y)$. We show that in any case the point $\xi$ belongs to the closure of the horosphere $\mathcal{HS}(\Phi, y)$ in $\partial_gX$ as well.

Let $c: \R_+ \to X$ be the parameterization of the ray $[y \xi]$. Let $x \in X$ be a point, such that $\Phi(x) > \Phi(y)$. For any $t > 0$ there exists a point $z_t$ in the horosphere $\mathcal{HS}(\Phi, x)$ with
$$|z_t c(t)| = \Phi(z_t) - \Phi(c(t)) = \Phi(x) -\Phi(c(t)).$$
This point $z_t$ belongs to the intersection of the horosphere $\mathcal{HS}(\Phi, x)$ with arbitrary ray from $c(t)$ in the direction opposite to the ray $[c(t) \pi_{hg}([\Phi])]$. The value $|y z_t| \ge t - |z_tc(t)|$ is non decreasing when $t \to \infty$, but the difference $\Phi(x) -\Phi(c(t))$ increases sublinearly: for any $k > 0$ and $b \in \R$ there exists $T = T(k, b) > 0$ such that 
$$\Phi(x) - \Phi(c(t)) < kt + b$$
for all $t > T$. Therefore 
$$\lim\limits_{t \to \infty} \frac{\Phi(x) - \Phi(c(t))}{t} = 0.$$
Any ray with start segment $[y z_t]$ passes out of the horoball $\mathcal{HB}(\Phi, y)$ and its endpoint lays out of $\mathcal{HB}_\infty([\Phi])$. At the same time, for any cone neighbourhood $\mathcal U$ of the point $\xi$ there exists sufficiently large $T > 0$, such that endpoints of rays with beginning part $[y z_t]$ belong to $\mathcal U$ for all $t > T$. This means that $\mathcal U$ has non empty intersection with the complement to the closure of the horoball $\mathcal{HB}(\Phi, y)$ in $\partial_gX$. This proves the inclusion $\xi \in \overline{\mathcal{HS}(\Phi, y)}_g$. 
\end{enumerate}
\qed

\begin{remark}\label{euklidLob}
The inverse inclusion to \eqref{horosphere_infty} can be false. The simplest example is Lobachevskii space, where all horospheres at infinity are empty but every horosphere of the space has an accumulation point at infinity --- its center. The geodesic completeness condition is essential here: it is easy to construct the situation when all level sets for a horofunction in non geodesically complete space are bounded, but the horosphere at infinity is not empty. For example, consider the horofunction $\Phi(x, y) = y$ on the part $y \ge \sqrt{|x|}$ of the  Euclidean plane with coordinates $(x, y)$
\end{remark}

\begin{remark}
It follows from the Lemma \ref{hb=hb} that all horoballs as sublevel sets for the horofunction $\Phi$ have the same boundary at infinity 
$$\partial_\infty \left(\mathcal{HB}(\Phi, x))\right)  = \overline{\mathcal{HB}(\Phi, x)}_g \setminus \mathcal{HS}(\Phi, x).$$  
Consequently, the horoball at infinity $\mathcal{HB}_\infty([\Phi])$ can be defined by the equality \eqref{hotoball_infty}. In another words, the horoball at infinity $\mathcal{HB}_\infty(\Phi)$ is the inverse limit for the system of closures $\overline{\mathcal{HB}(\Phi, t)}_g$ with inclusions $\overline{\mathcal{HB}(\Phi, t_1)}_g \subset \overline{\mathcal{HB}(\Phi, t_2)}_g$ when $t_1 \le t_2$ under $t \to -\infty$.
\end{remark}

The statement of the following theorem is another formulation of the Lemma \ref{limmm} in terms of horoballs at infinity.

\begin{theorema}\label{horlimits}
Let the consequence of horofunctions $\{\Phi_n\}_{n = 1}^\infty$ converges in compact-open topology to the horofunction $\Phi$ and the point $\xi \in \partial_gX$ is the limit of the consequence $\{\xi_n\}_{n = 1}^\infty \subset \partial_gX$ in the sense of the cone topology in $\partial_gX$, where
\begin{equation}\label{seqin}
\xi_n \in \mathcal{HB}_\infty(\Phi_n).
\end{equation}
Then $\xi \in \mathcal{HB}_\infty(\Phi)$.
\end{theorema}

\begin{remark}
The claim of the Theorem \ref{horlimits} is equivalent the inclusion
\begin{equation}\label{incllims}
\lim_{n \to \infty} \mathcal{HB}_\infty(\Phi_n) \subset \mathcal{HB}_\infty\left(\lim_{n \to \infty} \Phi_n\right),
\end{equation}
under the condition that horofunctions $\Phi_n$ converge to the horofunction $\Phi$.
Here the limit $\lim_{n \to \infty} \mathcal{HB}_\infty(\Phi_n)$ is the union of accumulation points for all different sequences $\{\xi_n\}_{n=1}^\infty$, where $\xi_n \in \mathcal{HB}_\infty(\Phi_n)$, converging in the sense of the cone topology on $\partial_gX$.
The example described in the Remark \ref{euklidLob} shows that the inclusion \eqref{incllims} can be strict.
\end{remark}

\textbf{Acknowledgement.} The author is very grateful to Sergey Vladimirovich Buyalo for his attention to the first version of the paper and the number of important remarks and corrections.

\noindent 

{\it Pomor State University,}
Arkhangelsk, Russia\newline
{\it Email:} \quad  {\tt  pdandreev@mail.ru}

\end{document}